\documentclass[12pt,oneside,english]{amsart}
\usepackage[latin1]{inputenc}
\usepackage{amssymb}

\makeatletter

\numberwithin{equation}{section}
\usepackage{babel}

\makeatother
\begin{document}
\baselineskip=17pt
\title[ Banach problem and a congruence for primes]{Banach matchboxes problem and a congruence for primes}

\author{Vladimir Shevelev}
\address{Department of Mathematics \\Ben-Gurion University of the
 Negev\\Beer-Sheva 84105, Israel. e-mail:shevelev@bgu.ac.il}


\begin{abstract}
Using an identity arising in the known Banach probability problem on matchboxes, we prove an unexpected congruence for odd prime $p:$ for $1\leq k\leq \frac{p-1}{2},\enskip \sum_{i=1}^{p-2k-1}2^{i-1}\binom{k-1+i}{k}\equiv 0\pmod p.$
\end{abstract}

\maketitle

\section{Introduction}
A classic probability Banach problem is the following (see Feller \cite{1}, Ch. 6, Section 8).
``Suppose a mathematician carries two matchboxes at all times: one in his left pocket and one in his right. Each time he needs a match, he is equally likely to take it from either pocket. Suppose he reaches into his pocket and discovers that the box picked is empty. If it is assumed that each of the matchboxes originally contained $n$ matches, what is the probability that there are exactly k matches in the other box?" \newline
The required probability is
\begin{equation}\label{1.1}
u_n(r)=\binom {2n-r} {n}2^{r-2n}, \enskip r=0,1,...,n.
\end{equation}
Since $\sum_{r=0}^n u_n(r)=1,$ then from (\ref{1.1}) we find
\begin{equation}\label{1.2}
\sum_{r=0}^n \binom {2n-r} {n}2^{r-n}=2^n.
\end{equation}
Riordan (\cite{2},\enskip Ch.1, problem 7) gave an independent proof of (\ref{1.2}).\newline
Using (\ref{1.2}), in this paper we prove an unexpected congruence for odd prime $p:$
\begin{equation}\label{1.3}
\sum_{i=1}^{p-2k-1}2^{i-1}\binom{k-1+i}{k}\equiv 0\pmod p,  \enskip  k=1,2,...,\frac{p-1}{2}.
\end{equation}

\section{Proof of the congruence}
Let us write (\ref{1.3}) in the form
$$\sum_{i=1}^{p-2k-1}2^{i-1}i(i+1)...(i+k-1)=(\sum_{i=1}^{p-2k-1}x^{i+k-1})^{(k)}\mid_{x=2}=$$
$$((x^{p-1+k}-x^k)(x-1)^{-1})^{(k)}\mid_{x=2}=$$ $$(\sum_{j=0}^k \binom{k}{j}(x^{p-1+k}-x^k)^{(j)}((x-1)^{-1})^{(k-j)})\mid_{x=2}. $$
We should prove that
$$\sum_{j=0}^k(-1)^{k-j}\binom{k}{j}(k-j)!(p-1-k)(p-2-k)...(p-k-j)2^{p-1-k-j}\equiv $$ $$\sum_{j=0}^k(-1)^{k-j}\binom{k}{j}(k-j)!k(k-1)...(k-j+1)2^{k-j}\pmod p.$$
Since $2^{p-1}\equiv1\pmod p,$ then this congruence is reduced to the identity
$$\sum_{j=0}^k\binom{k}{j}(k-j)!(k+1)(k+2)...(k+j)2^{-k-j}=$$ $$k!\sum_{j=0}^k(-1)^{j}\binom{k}{j}2^{k-j}=k!,$$
or
$$\sum_{j=0}^k 2^{-j}\binom{k}{j}(k-j)!(k+j)!=2^k k!^2,$$
or
$$\sum_{j=0}^k 2^{-j}\binom{k+j}{k}=2^k,$$
or, finally, putting here $j=k-i,$ we reduce congruence (\ref{1.3}) to the Banach identity (\ref{1.2})
$$\sum_{i=0}^k \binom{2k-i}{k}2^{i-k}=2^k.$$
This completes the proof.

\end{document}